\documentclass[a4paper, 12pt]{article}

\usepackage[koi8-r,cp1251]{inputenc}
\usepackage[english,russian]{babel}
\usepackage{amsfonts,amssymb,amsmath, hyperref}
\usepackage[final]{epsfig}
\usepackage{graphicx}

\newtheorem{theorem}{Теорема}
\newtheorem{lemma}{Лемма}
\newtheorem{definition}{Определение}
\newtheorem{corollary}{Следствие}
\newtheorem{remark}{Замечание}

\begin{document}

\begin{center}
\Large{Свойства функций Бернштейна нескольких комплексных переменных}
\end{center}

\

\begin{center}
А.\,Р.~Миротин
\end{center}

\begin{center}
amirotin@yandex.ru
\end{center}

\

Вводится многомерное обобщение класса функций Бернштейна и
исследуются свойства функций введенного класса. В частности,
дается новое доказательство интегрального представления функций
Бернштейна многих переменных. Рассматриваются примеры.

Ключевые слова: функция Бернштейна, абсолютно монотонная функция,
марковский процесс.

\

Abstract.
A multidimensional generalization of the Bernstein class of functions and the
properties of functions of the  introduced class are examined. In particular, a
new proof of the integral representation of Bernstein functions
of many variables is given. Examples are considered.

\


\

Mathematical Notes, 2013, Vol. 93, No. 2, pp. 257–265.
Original Russian Text © A. R. Mirotin, 2013, published in Matematicheskie Zametki, 2013, Vol. 93, No. 2, pp. 216–226.

\subsection{Введение}
\label{subsec1}

Функции Бернштейна одного переменного возникают в
    математическом анализе, теории потенциала
    (см. \cite{BCR}, с. 141 -- 142),  теории вероятностей (см.
    \cite{Fel}, глава XIII, \S 7; \cite{App}),
    а также теории операторов (см.,  например,  \cite{BBD}).
    Класс функций Бернштейна нескольких  переменных был
введен автором в \cite{Mir1} и использовался в работах \cite{Mir1}
--- \cite{Mir3} для
построения функционального исчисления генераторов
многопараметрических полугрупп операторов. Эти функции
 естественным образом возникают также в теории
марковских процессов. В заметке устанавливается ряд свойств
функций этого класса, включая их интегральное представление, даны
новые способы построения функций Бернштейна одного и двух
переменных, а также указана связь этого класса с теорией
марковских процессов.

Следует отметить, что в литературе чаще всего  встречается класс
~$\mathcal{B}$ (положительных) функций Бернштейна одного
переменного, теория которого тождественна теории введенного ниже
класса~$\mathcal {T}_1$ отрицательных функций Бернштейна,
поскольку элементы последнего класса получаются из функций
класса~$\mathcal{B}$ отражением относительно начала координат.

Ниже все меры считаются мерами Радона, неравенства между векторами
-- покоординатные. Для векторов $a,b\in \mathbb{C}^n$ выражение
$a\cdot b$ обозначает $\sum_{i=1}^n a_ib_i$.

\subsection{Основные свойства функций класса~${\cal T}_n$}
\label{subsec2}

\begin{definition}
 \label{definition1} \cite{Mir1} Скажем, что
неположительная функция  из $C^\infty((-\infty;0)^n)$ {\it
принадлежит классу ${\cal T}_n$}, если все ее частные производные
первого порядка абсолютно монотонны (функция из
$C^{\infty}((-\infty;0)^n)$ называется {\it абсолютно монотонной,}
если она неотрицательна вместе со своими частными производными
всех порядков).
\end{definition}

Последнее условие на функцию $\psi\in{\cal T}_n$ равносильно тому,
что $\partial^\alpha\psi\geq 0$ для любого мультииндекса
$\alpha\ne 0$. Функции класса ${\cal T}_n$ будем называть также
{\it отрицательными функциями Бернштейна $n$ переменных}.

Поскольку функция $\psi\in{\cal T}_n$  монотонно возрастает по
каждому переменному в отдельности, неравенство $x\leq y$, где
$x,y\in (-\infty;0)^n$, влечет $\psi(x)\leq \psi(y)$. Очевидно,
что ${\cal T}_n$ есть конус относительно поточечного сложения
функций и умножения на скаляр. Он также инвариантен относительно
сдвигов на векторы из $(0;\infty)^n$.  Более того, в силу
отмеченной выше  монотонности функция $\psi(s-c)-\psi(-c)$
принадлежит ${\cal T}_n$ вместе с $\psi\ (c\in (0;\infty)^n)$.

Следующая теорема  была анонсирована в \cite{Mir1} и доказана в
\cite{Mir3} с использованием одного нетривиального результата из
\cite{BCR} (относительно одномерного случая см.,  например,
\cite{Ah}, с. 256--258). Ниже мы даем прямое доказательство.
Остальные утверждения данной работы являются следствиями этой
теоремы {\rm(}далее запись $s\to -0$ означает, что $s_1\to -0,
\ldots, s_n\to -0${\rm)}.

\begin{theorem}
\label{theorem1}
 Пусть $m, n$ -- натуральные числа. Для
функции $\psi:(-\infty;0)^n\to (-\infty;0]$ следующие утверждения
равносильны:

{\rm 1)} $\psi\in {\cal T}_n;$

{\rm 2)} при $s\in(-\infty;0)^n$
\begin{equation}
\psi(s)=c_0+c_1\cdot s+\int\limits_{\Bbb{R}_+^n\setminus \{0\}}
\left(e^{s\cdot u}-1\right)d\mu (u), \label{2.1}
\end{equation}
  где
$$
c_0=\psi(-0):=\lim\limits_{s\to -0}\psi(s),\quad
c_1=\nabla\psi(-\infty)\in \Bbb{R}_+^n,
$$
а положительная мера
$\mu$ на $\Bbb{R}_+^n\setminus \{0\}$ определяется однозначно и обладает следующими
свойствами: для достаточно малого $\delta >0$ функции
$u\mapsto u_j\  \mu$-интегрируемы на
$[0;\delta)^n\setminus\{0\}, j=1,\ldots,n,$
и $\mu(\Bbb{R}_+^n\setminus[0;\delta)^n)<\infty$;

{\rm 3)} для любой  функции $f$, абсолютно монотонной на
$(-\infty;0)^m$, сложная функция $(s,r)\mapsto f(\psi(s),r)$
абсолютно монотонна на $(-\infty;0)^{m+n-1}\ (s\in
(-\infty;0)^n,r\in(-\infty;0)^{m-1}).$
\end{theorem}

\medskip
Доказательство. Утверждение $1\Rightarrow 3$ проверяется
непосредственным дифференцированием.

Для доказательства импликации $3 \Rightarrow 1$ достаточно
рассмотреть случай $f(x)=e^{vx_1}, v>0$.

Утверждение $2\Rightarrow 1$ доказывается дифференцированием под
знаком интеграла в (2.1) (условия, накладываемые на меру, гарантируют его сходимость).

Докажем, что  $1\Rightarrow 2$. В силу многомерного варианта
теоремы Бернштейна-Уиддера  абсолютно монотонная функция
$\partial\psi/\partial s_i$ есть $n$-мерное преобразование Лапласа
положительной меры $\sigma_i$, т.е. имеет интегральное
представление $(i=1,\ldots,n)$

$$
\frac{\partial\psi(s)}{\partial
s_i}=\int\limits_{\Bbb{R}_+^n}e^{s\cdot u}d\sigma_i(u).
$$

 Введем в рассмотрение множества
$U_i=\{u\in\Bbb{R}_+^n|u_i>0\}$ и определим в $U_i$ меру Радона
$\mu_i$ по правилу
$$
\int\limits_{U_i}\varphi
d\mu_i=\int\limits_{U_i}\varphi(u)\frac{1}{u_i}d\sigma_i(u),
$$

\noindent где $\varphi$ есть непрерывная функция с компактным
носителем, содержащимся в $U_i$. Заметим, что $U_i$ образуют
открытое покрытие пространства $\Bbb{R}_+^n\setminus \{0\}$. Кроме
того, меры $\mu_i$ согласованы в том смысле, что для любой пары
$(i,j)$ имеет место равенство сужений $\mu_i|(U_i\cap
U_j)=\mu_j|(U_i\cap U_j)$, то есть

$$
\frac{1}{u_i}d\sigma_i(u)=\frac{1}{u_j}d\sigma_j(u)\  { \mbox в }
\  U_i\cap U_j.
$$
\noindent В самом деле, $u_id\sigma_j(u)=u_jd\sigma_i(u)$,
поскольку обе части являются представляющими мерами для смешанной
производной $\partial^2\psi/\partial s_i\partial s_j
=\partial^2\psi/\partial s_j\partial s_i$. Следовательно (см.,
например, \cite{BurbInt}, глава 3, \S  2, Предложение 1),
существует такая положительная мера $\mu$ на $\Bbb{R}_+^n\setminus
\{0\}$, что $\mu|U_i=\mu_i$.

Предположим на время, что $\partial\psi(-0)/\partial s_i\ne\infty$
при всех $i$. Покажем, что тогда функция $u\mapsto e^{s\cdot u}-1$
$\mu$-интегрируема для любого $s\in (-\infty;0)^n$. Для $\delta>0$
рассмотрим множества $E_i=\{u\in\Bbb{R}_+^n|u_i\geq\delta\}$. Так
как при всех $i$ справедливо неравенство

$$
\int\limits_{E_i}e^{s\cdot u} d\mu=\int\limits_{E_i}e^{s\cdot u}
\frac{1}{u_i}d\sigma_i(u)
\leq\frac{1}{\delta}\int\limits_{E_i}e^{s\cdot u}
d\sigma_i(u)\leq\frac{1}{\delta} \frac{\partial\psi(s)}{\partial
s_i},
$$
\noindent то функция $u\mapsto e^{s\cdot u}$ $\mu$-интегрируема на
множестве $\cup_{j=1}^nE_j= \Bbb{R}_+^n\setminus [0;\delta)^n$.
Полагая в этом неравенстве последовательно $s_1\to-0,\ldots,
s_n\to-0$, получим
$\int_{E_j}d\mu\leq\delta^{-1}\partial\psi(-0)/\partial
s_i\ne\infty$ при всех $i$, а потому и функция $u\mapsto e^{s\cdot
u}-1$ $\mu$-интегрируема на $\Bbb{R}_+^n\setminus [0;\delta)^n$.
Далее, функция $u\mapsto (-s)\cdot u\ \mu$-интегрируема на
$[0;\delta)^n\setminus \{0\}$, поскольку меры
$u_jd\mu(u)=d\sigma_j(u)$ конечны на этом множестве. Тогда
соотношение $1-e^{s\cdot u} \sim (-s)\cdot u\ (u\to 0)$
показывает, что функция $u\mapsto e^{s\cdot u}-1$
$\mu$-интегрируема на множестве $[0;\delta)^n\setminus \{0\}$ при
достаточно малом $\delta>0$. Окончательно получаем, что эта
функция $\mu$-интегрируема на множестве
$\Bbb{R}_+^n\setminus\{0\}$.

Следовательно, при всех $i$

$$
\frac{\partial\psi(s)}{\partial
s_i}=\int\limits_{\Bbb{R}_+^n}e^{s\cdot u}d\sigma_i(u)=
\frac{\partial}{\partial
s_i}\left(\int\limits_{\Bbb{R}_+^n\setminus \{0\}}\left(e^{s\cdot
u}-1\right)) d\mu(u)+c_1^is_i\right),
$$
\noindent где
$c_1^i=\sigma_i(\{0\})=\partial\psi(-\infty)/\partial s_i$, откуда
и следует (2.1).

Наконец, если $\partial\psi(-0)/\partial s_i=\infty$ при некотором
$i$, то для $\epsilon>0$ рассмотрим вектор
$\overline{\epsilon}:=(\epsilon,\ldots,\epsilon)$ и функцию
 $\psi_\epsilon(s):=\psi(s-\overline{\epsilon})$ из ${\cal T}_n$, для которой
$\partial\psi_\epsilon(-0)/\partial s_i\ne\infty$ при всех $i$.
Поскольку

$$
\frac{\partial\psi_\epsilon(s)}{\partial
s_i}=\int\limits_{\Bbb{R}_+^n} e^{(s-\overline{\epsilon})\cdot
u}d\sigma_i(u)=\int\limits_{\Bbb{R}_+^n} e^{s\cdot
u}e^{-\overline{\epsilon}\cdot u}d\sigma_i(u)
$$

\noindent есть преобразование Лапласа меры
$d\sigma_i^{\epsilon}(u) =e^{-\overline{\epsilon}\cdot
u}d\sigma_i(u)$,
 и $\psi_\epsilon(0)=\psi(-\overline{\epsilon})$, $\sigma_i^{\epsilon}(\{0\})=
\sigma_i(\{0\})$,  интегральное представление (2.1) для
$\psi_\epsilon$ имеет вид

$$
\psi(s-\overline{\epsilon})=\psi(-\overline{\epsilon})+c_1\cdot
s+\int\limits_{\Bbb{R}_+^n\setminus \{0\}} \left(e^{s\cdot
u}-1\right)e^{-\overline{\epsilon}\cdot u}d\mu (u).
$$
\noindent Для завершения доказательства формулы (2.1) осталось
положить тут $\epsilon\to +0$ и сослаться на теорему Б. Леви. Требуемые свойства меры $\mu$
установлены в  лемме 3.1 из \cite{Mir3}.

Наконец, единственность меры $\mu$ следует из теоремы
единственности для преобразования Лапласа, поскольку
дифференцирование (2.1) под знаком интеграла показывает, что
$\partial\psi/\partial s_i$ есть преобразование Лапласа меры
$u_id\mu(u)$.

\bigskip
Композиция функций сохраняет принадлежность к классу ${\cal T}:=
\cup_{n=1}^\infty{\cal T}_n.$ Точнее, имеет место такое следствие,
показывающее, что ${\cal T}$ есть операда в смысле \cite{MSS}.

\begin{corollary}
 Если $\psi_1(s)\in {\cal T}_n$ и $\psi_2(r)\in {\cal
T}_m,$ то функция $\psi(r,'s):=\psi_1(\psi_2(r),'s)$ принадлежит
${\cal T}_{m+n-1}$ (мы полагаем $'s=(s_2,\ldots,s_n)).$ В
частности, класс ${\cal T}_1$ устойчив относительно композиции.
\end{corollary}

Доказательство. По теореме 1 для любой абсолютно монотонной на
$(-\infty;0)$ функции $f$ функция $f_1(s):= f(\psi_1(s))$
абсолютно монотонна. Аналогично, функция
$f(\psi(r,'s))=f_1(\psi_2(r),'s) $ также абсолютно монотонна, и
снова в силу теоремы 1 функция $\psi(r,'s)$ принадлежит классу
${\cal T}_{m+n-1}$.

\bigskip
\begin{remark}
 Каждая функция $\psi\in{\cal T}_n$ по формуле

\begin{equation}
\psi(z)=c_0+c_1\cdot z+\int\limits_{\Bbb{R}_+^n\setminus \{0\}}
(e^{z\cdot u}-1)d\mu (u)     \label{2.2}
\end{equation}
продолжается до функции, голоморфной в области
$$\{{\rm
Re}z<0\}:=\{z\in\mathbb{C}^n:{\rm Re}z_j<0,j=1, \ldots,n\}
$$
и непрерывной в ее замыкании
(абсолютная сходимость интеграла в (2.2) была доказана в
\cite{Mir3}). Класс функций, возникающий в результате такого аналитического продолжения,
будем по-прежнему обозначать ${\cal T}_n$.
\end{remark}

\bigskip
Далее $S_\theta$ есть замыкание сектора
$$
\left\{\zeta\in\Bbb{C}:\frac{\pi}{2}+\frac{\theta}{2} <{\rm
arg}\zeta<\frac{3\pi}{2}-\frac{\theta}{2}\right\},\quad\theta\in[0;\pi].
$$

\begin{theorem}
\label{theorem2}Любая  функция из ${\cal T}_n$, отличная от
тождественно нулевой, отображает область $\{{\rm Re}z<0\}$  в
полуплоскость $\{{\rm Re}\zeta<0\}$,  а конус $S_\theta^n$ -- в
сектор $S_\theta$.
\end{theorem}

Доказательство. Полагая в (2.2) $z=s+iy$, получаем, что
\begin{equation}
{\rm Re}\psi(z)=c_0+c_1\cdot
s+\int\limits_{\Bbb{R}^n_+\setminus\{0\}} (e^{s\cdot u}\cos(y\cdot
u)-1)d\mu(u).\label{2.3}
\end{equation}

Поскольку каждое слагаемое в правой части (2.3) неположительно, то
${\rm Re}\psi(z)\leq 0$ при всех $z$ c ${\rm Re}z<0$, причем
равенство достигается лишь когда каждое слагаемое в правой части
(2.3) равно нулю. Интерес представляет случай $\mu\ne 0$. В этом
случае из равенства ${\rm Re}\psi(z)=0$ следует, что
$c_0=c_1=0,\quad e^{s\cdot u_0}\cos(y\cdot u_0)=1$ при некотором
$u_0\in \Bbb{R}_+^n\setminus\{0\}$. Последнее равенство
невозможно, поскольку $s<0$, что и доказывает первое утверждение
теоремы.

Для доказательства второго утверждения возьмем $z=s+iy\in
S_\theta^n,z\ne 0,$ т.е. $s_j<0$ и $|y_j|/(-s_j)\leq {\rm
ctg}(\theta/2),j=1,\ldots,n.$ Тогда

$$
\frac{|e^{u\cdot s}\sin(u\cdot y)|}{1-e^{u\cdot s}\cos(u\cdot y)}=
\frac{|\sin(u\cdot y)|}{e^{-u\cdot s}-\cos(u\cdot y)}\leq
\frac{|u\cdot y|}{(1-u\cdot s)-1}\leq {\rm ctg}\frac{\theta}{2}.
$$

Теперь равенство (2.3) вместе с равенством

$$
{\rm Im}\psi(z)=c_1\cdot y+\int\limits_{\Bbb{R}_+^n\setminus\{0\}}
e^{s\cdot u}\sin(y\cdot u)d\mu(u),
$$

\noindent которое также следует из (2.2), показывают, что $|{\rm
Im}\psi(z)|\leq {\rm ctg}(\theta/2) (-{\rm Re}\psi(z)),$ т.е.
$\psi(z)\in S_\theta$, что и требовалось доказать.

\begin{definition}
 \label{definition2}Будем говорить, что последовательность
$\psi_k$ функций из ${\cal T}_n$ {\it сходится  в ${\cal T}_n$} к
некоторой функции $\psi$, если $\psi$ определена на множестве
$(-\infty;0)^n$ и $\psi_k$ сходятся к $\psi$ поточечно на этом
множестве.
\end{definition}

Следующий результат показывает, что описанная сходимость для
класса ${\cal T}_n$ является естественной.

\begin{theorem}
\label{theorem3}Пусть последовательность $\psi_k$ сходится в
${\cal T}_n$ к функции $\psi$. Тогда

{\rm 1)} $\psi\in{\cal T}_n;$

{\rm 2)} последовательность $\psi_k$ сходится к $\psi$ равномерно
на компактных подмножествах области $\{{\rm Re}z<0\}$;

{\rm 3)}  последовательность $\psi_k$ равномерно ограничена на
компактных подмножествах конуса $S_\theta^n, \theta\in (0;\pi)$.
\end{theorem}

Доказательство. Как отмечено в доказательстве теоремы 2, ${\rm
Re}\psi_k(z) \leq 0$ при $k=1,2,\ldots, {\rm Re}z\leq 0$, т. е.
$\psi_k$ выпускают значения из полуплоскости $\{{\rm Re}\zeta
>0\}$. Следовательно, в силу многомерного варианта теоремы Монтеля
семейство $(\psi_k)$ нормально, т. е. любая его
подпоследовательность содержит подпоследовательность, сходящуюся
равномерно внутри области $\{{\rm Re}z<0\}$ к аналитической
функции или к бесконечности. Однако последнее невозможно ввиду
поточечной сходимости $\psi_k$. Пределы всех сходящихся равномерно
внутри $\{{\rm Re}z<0\}$ подпоследовательностей последовательности
$\psi_k$ совпадают с $\psi$ на $(-\infty;0)^n$, а потому по одной
из форм теоремы единственности совпадают между собой на $\{{\rm
Re}z<0\}$ (см., например, \cite{Sh}, c. 32). Но тогда и сама
последовательность $\psi_k$ сходится равномерно внутри $\{{\rm
Re}z<0\}$ к функции, являющейся аналитическим продолжением $\psi$.
Это продолжение мы также обозначим $\psi$. Теперь по теореме
Вейерштрасса для любого мультииндекса $\alpha$ при $k\to\infty$
имеет место сходимость
$\partial^\alpha\psi_k(s)\to\partial^\alpha\psi(s),s\in
(-\infty;0)^n$, а потому $\psi\in{\cal T}_n$, что доказывает
первые два утверждения  теоремы.

Для доказательства третьего утверждения нам понадобится

\begin{lemma}
\label{lemma1}Для любой функции $\psi \in{\cal T}_n$ с
интегральным представлением {\rm (2.2)} справедливо неравенство
\begin{equation}
{\rm Re}\psi(s+iy)-\psi(s)\geq 2({\rm Re}\psi(r+iy)-c_0-c_1\cdot
r),\label{2.4}
\end{equation}
 где $y\in\Bbb{R}^n;r,s\in (-\infty;0)^n.$
\end{lemma}

Доказательство леммы. В силу формулы (2.3)
$$
{\rm Re}\psi(s+iy)=c_0+c_1\cdot
s+\int\limits_{\Bbb{R}^n_+\setminus\{0\}} (e^{s\cdot u}\cos(y\cdot
u)-1)d\mu(u)=
$$
$$
c_0+c_1\cdot s+\int\limits_{\Bbb{R}^n_+\setminus\{0\}} (e^{s\cdot
u}-1)\cos(y\cdot u)d\mu(u)
+\int\limits_{\Bbb{R}^n_+\setminus\{0\}} (\cos(y\cdot
u)-1)d\mu(u)\leq 0.
$$
\noindent Следовательно,

$$
|{\rm Re}\psi(s+iy)|=|c_0|+|c_1\cdot
s|+\int\limits_{\Bbb{R}^n_+\setminus\{0\}} (1-e^{s\cdot
u})|\cos(y\cdot u)|d\mu(u)+\\
$$

$$
\int\limits_{\Bbb{R}^n_+\setminus\{0\}} (1-\cos(y\cdot
u))d\mu(u)\leq
$$
$$
-c_0-c_1\cdot s+\int\limits_{\Bbb{R}^n_+\setminus\{0\}}
(1-e^{s\cdot u})d\mu(u) +\int\limits_{\Bbb{R}^n_+\setminus\{0\}}
(1-\cos(y\cdot u))d\mu(u)=
$$
\begin{equation}
-\psi(s)+\int\limits_{\Bbb{R}^n_+\setminus\{0\}} (1-\cos(y\cdot
u))d\mu(u).\label{2.5}
\end{equation}
\noindent Так как все
слагаемые в правой части  и подинтегральная функция в (2.3)
неположительны, то
\begin{equation}
|{\rm Re}\psi(r+iy)|=-c_0-c_1\cdot
r+\int\limits_{\Bbb{R}^n_+\setminus\{0\}} (1-e^{r\cdot
u}\cos(y\cdot u))d\mu(u).\label{2.6}
\end{equation} \noindent

Заметим, что при $t\in (0;1],z\in[-1;1]$ справедливо неравенство
$1-z\leq 2(1-tz)$. Следовательно, $1-\cos(y\cdot u)\leq
2(1-e^{r\cdot u} \cos(y\cdot u))$, и в силу (2.6) имеем
$$
\int\limits_{\Bbb{R}^n_+\setminus\{0\}} (1-\cos(y\cdot
u))d\mu(u)\leq 2\int\limits_{\Bbb{R}^n_+\setminus\{0\}}
(1-e^{r\cdot u}\cos(y\cdot u))d\mu(u)=
$$
$$
2(|{\rm Re}\psi(r+iy)|+c_0+c_1\cdot r).
$$
\noindent С учетом этого неравенства (2.5) дает
$$
|{\rm Re}\psi(s+iy)|\leq -\psi(s)+2(|{\rm
Re}\psi(r+iy)|+c_0+c_1\cdot r),
$$
т.е.
$$
-{\rm Re}\psi(s+iy)\leq -\psi(s)+2(-{\rm
Re}\psi(r+iy)+c_0+c_1\cdot r),
$$
\noindent что и завершает доказательство леммы.

\medskip

Доказательство утверждения 3). Любое ограниченное множество $M$ из
$S_\theta^n$ при некотором $a<0$ содержится в множестве
$\Delta^n$, где $\Delta=\{\zeta\in S_\theta:{\rm Re}\zeta\geq
a\}$.  В силу неравенства (2.4) при $r=(a,\ldots,a)$ имеем
$$
|{\rm Re}\psi_k(s+iy)|\leq |\psi_k(s)|+2|{\rm
Re}\psi_k(r+iy)|+\nabla\psi_k(-\infty)\cdot r+2\psi_k(0)\leq
$$
\begin{equation}
|\psi_k(s)|+2|{\rm Re}\psi_k(r+iy)|+2\psi_k(0).\label{2.7}
\end{equation}
 Далее, последовательность $\psi_k(a,\ldots,a)$, будучи сходящейся,
ограничена. А так как функция $s\mapsto \psi_k(s_1,\ldots,s_n)$
возрастает по каждому переменному в отдельности, то для всех
$s\geq (a,\ldots,a)$ выполняется неравенство
$\psi_k(a,\ldots,a)\leq\psi_k(s)\leq 0,$ т.е. последовательность
$\psi_k$ равномерно ограничена на $[a;0]^n$. Из утверждения 2)
следует, что последовательность $\psi_k$ равномерно ограничена на
компакте $K=\{z\in \Delta^n:{\rm Re}z_j=a,j=1,\ldots,n\}$. По
доказанному выше первое и третье слагаемые в правой части (2.7)
 ограничены, когда $s+iy\in\Delta^n$, так как тогда
$s\in [a;0]^n$. Второе же слагаемое равномерно ограничено,
поскольку $r+iy\in K$.  Теперь утверждение 3) сразу следует из
неравенства $|{\rm Im}\psi(z)|\leq {\rm ctg}(\theta/2) (-{\rm
Re}\psi(z))$, установленного при доказательстве теоремы 2.

\subsection{Примеры}
\label{subsec3}

Функции $c^\alpha-(c-s)^\alpha\ (0<\alpha<1, c\geq 0),\ \log
b-\log(b-s),\ {\rm arch b}-{\rm arch}(b-s)\ (b\geq 1)$ принадлежат
классу ${\cal T}_1$ (см., например, \cite{Mir3} для случаев $c=0,
b=1$; общий случай следует из замечания после определения 1).
Приведем теорему, которая позволяет строить новые примеры функций
из этого класса.

\begin{theorem}
 \label{theorem4}Если функция $w(-p)$ есть изображение по
Лапласу неотрицательного неубывающего оригинала, а $\psi\in {\cal
T}_1$, причем $\psi'(-\infty)=0$, то и функция
\begin{equation}
\varphi(s)=\int\limits_s^0\psi(t)w(t)dt,\quad s<0\label{3.1}
\end{equation}
также принадлежит ${\cal T}_1$ (при условии, что интеграл
существует).
\end{theorem}

Доказательство. Ясно, что $\varphi(s)\leq 0 \mbox{ при } s<0$, так
как $\psi(t)\leq 0, w(t)\geq 0 (t<0)$.  Рассмотрим два случая.

1) Пусть $\psi(s)=e^{us}-1, u\geq 0.$ Соответствующую ей по
формуле (3.1) функцию обозначим $\varphi_u$. Проверим, что функция
$\varphi_u'(s)=-\psi(s)w(s)$ абсолютно монотонна. Если
$w(-p)=({\cal L}f)(p),$ где $f\geq 0$ -- неубывающий оригинал,
(здесь и ниже ${\cal L}$ -- преобразование Лапласа), то
$$
\varphi'_u(s)=-e^{us}w(s)+w(s)=\int\limits_0^\infty
e^{rs}f(r)dr-\int\limits_0^ \infty e^{(u+r)s}f(r)dr=
$$
$$
\int\limits_0^\infty e^{rs}f(r)dr-\int\limits_u^ \infty
e^{rs}f(r-u)dr=\int\limits_0^u e^{rs}f(r)dr+\int\limits_u^ \infty
e^{rs}(f(r)-f(r-u))dr=
$$
$$
\int\limits_0^\infty e^{rs}g(r)dr,
$$
\noindent где $g(r)=f(r)-f(r-u)\geq 0,r\in\Bbb{R}$, (у нас
$f(s)=0$ при $s<0$). Таким образом, $\varphi'_u(s)=({\cal L}g)(-s)
(s<0)$ есть абсолютно монотонная функция.

2)  Функция $\psi\in {\cal T}_1$ имеет интегральное представление
(2.1), в котором $c_1=0$. Тогда в силу теоремы Фубини
$$
\varphi(s)=c_0\int\limits_s^0
w(t)dt+\int\limits_0^\infty\left(\int\limits_s^0
(e^{ut}-1)w(t)dt\right)d\mu(u)=
$$
\begin{equation}
 c_0\int_s^0 w(t)dt+\int\limits_0^\infty
\varphi_u(s)d\mu(u)\label{3.2}
\end{equation}
(теорема Фубини применима, так как подинтегральная функция
сохраняет знак в области интегрирования).

Первое слагаемое в правой части (3.2) принадлежит ${\cal T}_1$,
поскольку оно неположительно, и функция
$(c_0\int_s^0w(t)dt)'=(-c_0)w(s)$ абсолютно монотонна.

Представляя интеграл как предел интегральных сумм, получаем, что
второе слагаемое в (3.2) есть предел в ${\cal T}_1$
последовательности функций вида ($\alpha_j>0$)
$$
\sigma(s)=\sum\limits_{j=1}^m\varphi_{u_j}(s)\alpha_j,
$$
\noindent которые принадлежат ${\cal T}_1$ по доказанному в случае
1. По теореме 3 получаем теперь, что и второе слагаемое в (3.2)
принадлежит ${\cal T}_1$, что и завершает доказательство.

\bigskip
 Пример 1. Полилогарифмическая функция
порядка $p\in \Bbb{N}$ имеет вид
$$
{\rm Li_p}(z)=\sum\limits_{n=1}^\infty \frac{z^n}{n^p},\quad|z|<1
$$
\noindent и может быть аналитически продолжена в полуплоскость
$\{{\rm Re}z<0\}$, причем справедливо рекуррентное соотношение
$$
{\rm Li_p}(s)=\int\limits_s^0 {\rm Li_{p-1}}(t)(-t^{-1})dt
$$
\noindent (см., например, \cite{Gonch}). Как легко проверить,
${\rm Li_1}(s)=-\log(1-s)\in{\cal T}_1$. Кроме того,
$-t^{-1}=({\cal L}\theta)(-t),$ где $\theta$ -- функция Хевисайда.
Следовательно, теорема 4 и отмеченное рекуррентное соотношение
показывают, что ${\rm Li_p}\in{\cal T}_1$, лишь только ${\rm
Li_{p-1}}\in{\cal T}_1$. В силу принципа индукции это означает,
что ${\rm Li_p}\in{\cal T}_1$ при всех $p$.

Функции из ${\cal T}_n$ можно строить, используя суммы функций из
${\cal T}_1$ (от различных аргументов), а также операцию
композиции. Укажем еще один способ построения функций из ${\cal
T}_2$, обобщающий первое утверждение теоремы 11.1 из \cite{Mir3}.

\begin{theorem}
 \label{theorem5} Пусть функция $\psi_1\in {\cal T}_1$ имеет
интегральное представление
$$
\psi_1(s)=c_0+\int\limits_0^\infty (e^{su}-1)d\mu_1(u),
$$
причем $\omega:=\int_0^\infty ud\mu_1(u)\ne\infty$. Тогда функция
$$
\psi(s_1,s_2):=\frac{\psi_1(s_1)-\psi_1(s_2)}{s_1-s_2}-\omega
$$
принадлежит ${\cal T}_2$ и имеет интегральное представление
$$
\psi(s_1,s_2)=\int\limits_0^\infty\int\limits_0^\infty
(e^{s_1u_1+s_2u_2}-1) d\mu(u_1,u_2),
$$
\noindent где $d\mu(u_1,u_2)$ --- образ меры  $d\mu_1(v)dw$ при
отображении
$$
u_1=\frac{v+w}{2},\quad u_2=\frac{v-w}{2}.
$$
\end{theorem}

 (Под значением функции $\psi$ при $s_1=s_2$ мы, как
обычно, понимаем ее предел при $s_2\to s_1$).

Доказательство.  Пусть $\Delta$ есть  угол в  $(v,w)$-плоскости,
ограниченный биссектрисами первого и четвертого координатных
углов. Рассмотрим двойной интеграл
$$
\frac{1}{2}\int\!\int_\Delta
(e^{s_1(v+w)/2+s_2(v-w)/2}-1)d\mu_1(v)dw
$$
\begin{equation}
=\frac{1}{2}\int\limits_0^\infty d\mu_1(v)\int\limits_{-v}^v
(e^{s_1(v+w)/2+s_2(v-w)/2}-1)dw. \label{3.3}
\end{equation}
В случае $s_1\ne s_2$, вычисляя внутренний интеграл, получаем для
 правой части (3.3) значение
$$
\frac{1}{s_1-s_2}\int\limits_0^\infty
(e^{s_1v}-e^{s_2v})d\mu_1(v)-\omega =\psi(s_1,s_2).
$$
В случае $s_1=s_2$ аналогичные вычисления дают для этой части
значение
$$
\int\limits_0^\infty
(e^{s_1v}-1)vd\mu_1(v)=\psi_1'(s_1)-\omega=\psi(s_1,s_1),
$$

 С другой стороны, сделав в интеграле, стоящем в левой
части (3.3), замену переменных по формулам $v=u_1+u_2,w=u_1-u_2$,
получим, что он равен
$$
  \int\limits_0^\infty\int\limits_0^\infty(e^{s_1u_1+s_2u_2}-1)d\mu(u_1,u_2),
$$

что и завершает доказательство.

 Пример 2. Как отмечено выше, функция $\psi_1(s)=\log
b-\log(b-s)\ (b\geq  1)$ принадлежит классу    ${\cal T}_1$.
Нетрудно подсчитать, что этом случае  $\omega=1/b$. Следовательно,
функция
$$
\frac{1}{s_1-s_2}\log\frac{b-s_2}{b-s_1}-\frac{1}{b}
$$
принадлежит ${\cal T}_2$.

 В заключение покажем, как функции класса ${\cal T}_n$
  возникают в теории марковских процессов (при
$n=1$ это было установлено Филлипсом \cite{HiF}, теорема 23.15.1;
см.  также \cite{App}, с. 1339--1340). Рассмотрим однородный во
времени и пространстве (см., например, \cite{HiF}, п. 23.13)
марковский процесс $\xi$ в фазовом пространстве $\Bbb{R}_+^n$.
Соответствующее ему семейство вероятностей перехода $P(t,x), t>0,
x\in \Bbb{R}_+^n$ порождает в силу уравнения Чепмена-Колмогорова
сверточную полугруппу вероятностных мер $P_t, t>0$ на
$\Bbb{R}_+^n$.  Пусть $g_t(s):=({\cal L}P_t)(-s)$ (${\cal L}$
обозначает $n$-мерное преобразование Лапласа, $s<0$) --
односторонняя характеристическая функция функции распределения
$P(t,\cdot)$. Ясно, что  $g_t$ абсолютно монотонна. Кроме того,
 для любого $m\in \Bbb{N}$
$$
g_{mt}(s)=\int\limits_{\Bbb{R}^n_+} e^{s\cdot u}dP_{mt}(u)=
\int\limits_{\Bbb{R}^n_+} e^{s\cdot u}dP_t^{\ast m}(u)=(g_t(s))^m,
$$
\noindent где $P_t^{\ast m}$ обозначает $m$-ю сверточную степень
меры $P_t$.
 Заменяя в последнем равенстве $t$ на $(k/m)t$, получаем
   $g_{rt}(s)=(g_t(s))^r$ для любого положительного
рационального $r=k/m \ (k,m\in \mathbb{N})$. (Если дополнительно
предположить, что отображение $t\mapsto P_t$ узко непрерывно, то
при любом $a>0$ справедливо равенство $g_{at}(s)=(g_t(s))^a$. В
самом деле, если последовательность положительных рациональных
чисел $r_k$ сходится к $a$, то в силу нашего условия непрерывности
$g_{at}(s)=\lim_{k\to\infty} g_{r_kt}(s)=(g_t(s))^a$ (см.,
например, \cite{BurbInt}, с. 529, предложение 14). В частности,
$g_t(s)=(g_1(s))^t$.) Далее, вполне монотонная функция $g_1(-r)=
({\cal L}P_1)(r)$, как было показано выше, удовлетворяет равенству
$g_1(-r)^{1/m}=g_{1/m}(-r), m\in \mathbb{N}, r\in (0;\infty)^n$ и
значит ее представляющая мера является безгранично делимым распределением
(т. е. $g_1(-r)$ безгранично делима в смысле \cite{UMZ}). Поэтому $g_1(-r)=e^{\psi(-r)}$, где
$\psi\in{\cal T}_n$ (в одномерном случае это хорошо известно, см.,
например, \cite{Fel}, глава XIII, \S 7; в общем случае это
следует, например, из теоремы 3 в \cite{UMZ}). По аналогии с
одномерным случаем
 функцию $\psi$ естественно называть  лапласовским показателем  процесса $\xi$. Таким
образом, \emph{лапласовский показатель  рассматриваемого
марковского процесса принадлежит} ${\cal T}_n$.

Следует отметить, что роль лапласовского показателя в теории
марковских процессов весьма велика. Например, рассуждая как в
\cite{HiF}, теорема 23.15.2, где рассмотрен случай $n=1$, можно
показать, что генератор марковской полугруппы, отвечающей $\xi$,
имеет вид $\psi(-\nabla)$, где значение лапласовского показателя
$\psi$ от набора операторов $-\nabla =(-\partial/\partial
s_1,\ldots,-\partial/\partial s_n)$ понимается в смысле
исчисления, построенного в  \cite{Mir3}.



\begin{thebibliography}{99}


\bibitem{BCR}
  C.\,~Berg,  J.\,P.\,R.~Christensen, P.\,~
Ressel,
 Harmonic analysis on semigroups,
 Grad. Texts in Math., \textbf{100}, Springer-Verlag,
 Berlin,
  1982


\bibitem{Fel}
В.\,~Феллер,
 Введение в теорию
вероятностей и ее приложения,
  Т. 2, Мир,
 М.,  1984

\bibitem{App}
D.\,~Applebaum, Levy Processes -- From
Probability to Finance and Quantum Groups,
 Notices of the
American Mathematical Society,
vol 51,
2004,
 1336--1347

\bibitem{BBD}
C.\,~Berg, K.\,~Boyadzhiev, R.\,~deLaubenfelse,
Generation of generators of holomorphic semigroups,
 J. Austral. Math. Soc. (Series A),
 vol 55,
 1993,
  246--269


\bibitem{Mir1}
А.\,Р.~Миротин,
   Действие функций
класса Шенберга ${\cal T}$ на конусе диссипативных элементов
банаховой алгебры,
Матем. заметки, т. 61,  1997, 630--633; English
transl., Math. Notes 61:4 (1997), 524–527.

\bibitem{Mir2}
 А.\,Р.~Миротин, Функции класса Шенберга
${\cal T}$ действуют в конусе диссипативных элементов банаховой
алгебры, II,
 Матем. заметки, т.64, 1998,  423--430; English transl., Math. Notes 64:3 (1998), 364–370.

\bibitem{OaM}
 Mirotin, A.R. Bernstein functions of several semigroup generators on Banach spaces under bounded perturbations. Operators and Matrices. \textbf{11}, 199--217   (2017)

\bibitem{OaMII}
 Mirotin, A.R. Bernstein functions of several semigroup generators on Banach spaces under bounded perturbations. II. Operators and Matrices. \textbf{12},  445 -- 463  (2018)

\bibitem{SMZ2011}
 А.\,Р.~Миротин. О некоторых свойствах многомерного функционального исчисления Бохнера-Филлипса. Сиб. матем. журнал. \textbf{52} (6),  1300 -- 1312  (2011); English transl.: Siberian Mathematical Journal. \textbf{52} (6),  1032--1041   (2011)

\bibitem{Mir09}
А.\,Р.~Миротин. О многомерном функциональном исчислении Бохнера-Филлипса. ПФМТ, \textbf{1} (1), 63–-66  (2009) (Mirotin, A.R. On multidimensional Bochner-Phillips functional calculus. Probl. Fiz. Mat.
Tekh. \textbf{1} (1), 63–-66  (2009))

\bibitem{Izv}
 А.\,Р.~Миротин. О совместных спектрах наборов неограниченных операторов. Изв. РАН. Серия математическая.
\textbf{79} (6), 145 -- 170 (2015); English transl.: Izv.Math. \textbf{79} (6), 1235 -- 1259
(2015). DOI 10.1070/IM2015v079n06ABEH002779




\bibitem{Mir3}
  А.\,Р.~Миротин,
  Многомерное ${\cal T}$-исчисление от генераторов
$C_0$-полугрупп,
 Алгебра и анализ, т. 11, 1999, 142--170; English transl.: St. Petersburg Math. J. \textbf{11} (2), 315–-335  (1999)

\bibitem{Ah}  Н.\,И.~Ахиезер,  Классическая проблема
моментов и некоторые вопросы анализа, связанные с нею,
Физматгиз,  М., 1961




\bibitem{BurbInt}
 Н.~Бурбаки,  Интегрирование. Меры на
локально компактных пространствах. Продолжение меры.
Интегрирование мер. Меры на отделимых пространствах.,  Элементы  математики, Наука, М., 1977



\bibitem{MSS}
M.~Markl,  S.~Shnider , and J.~Stasheff,
Operads in Algebra, Topology and Physics  Math. Surveys Monogr., 96, Amer.
Math. Soc.,
 Providence, RI,
 2002


\bibitem{Sh} Б.~В.\,~Шабат  Введение в комплексный анализ. Ч. 2. Функции
нескольких переменных.
 Наука,
 М.  1976




\bibitem{Gonch}
  A.\,B.~Goncharov. The clasical
trilogarithm, algebraic K-theory of fields, and Dedekind zeta
function, Bull.  Amer. Math. Soc. (N.S.) vol.
24, 1991, 155--162

\bibitem{HiF}   Э.~Хилле,  Р.~Филлипс,  Функциональный
анализ и полугруппы,  ИЛ, M.,  1962



\bibitem{UMZ}
   А.\,Р.~Миротин,
  Вполне монотонные функции на полугруппах Ли,
  Укр. матем. журн., т.52,  2000,
  841--845; English transl.: Ukranian Math. J.,   vol 52,  2000.

\end{thebibliography}
\end{document}